\documentclass[11pt]{article}

\usepackage{wrapfig}
\usepackage{hyperref}
\usepackage[margin=1.0in]{geometry}
\usepackage{todonotes}

\usepackage{graphicx}
\usepackage{epstopdf}
\usepackage{amssymb,amsmath,amscd,amsthm,amsfonts}
\usepackage{graphicx,psfrag,epsfig,multirow}
\usepackage{cancel}
\usepackage{url}
\usepackage{caption}

\newtheorem{thm}{Theorem}

\newtheorem{rem}{Remark} \def\illustration #1 by #2 (#3) (#4){\leavevmode
  \vbox to #2{
   \hrule width #1 height 0pt depth 0pt
   \vfill
   \special{illustration #3 scaled #4}}}

\newcommand{\eps}{\varepsilon}
\newcommand{\p}{\partial}

\newcommand{\ber}{\mathbf{e}_r}
\newcommand{\bet}{\mathbf{e}_{_\theta}}

\begin{document} 
\title{The ponderomotive Lorentz force} 
\author{Graham Cox\footnote{Department of Mathematics and Statistics, Memorial University of Newfoundland. {\tt gcox@mun.ca}} \ and Mark Levi\footnote{Department of Mathematics, Pennsylvania State University. {\tt levi@math.psu.edu}}
} 
\bibliographystyle{plain} 

\maketitle 

\begin{abstract}
This paper describes a curious phenomenon: a particle in a rapidly varying potential is subject to an effective magnetic-like force. This force is in addition to the well-known ponderomotive force, but it has not been shown to exist before except for the linear case of a rapidly rotating quadratic saddle potential. We show that this is a universal phenomenon: the magnetic-like force arises generically in potential force fields  with rapid periodic time dependence, including but not limited to rotational dependence. 
\end{abstract}

\section{Introduction}  

\subsection*{The main point of the paper} 
\begin{figure}[thb]
 	\captionsetup{format=hang}
	\center{  \includegraphics{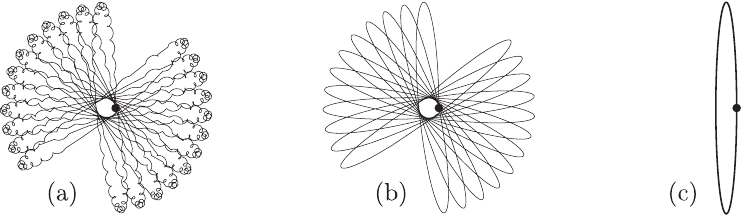}}
	\caption{Motion in a rotating quadratic saddle potential with $U_0(x)= (x_1 ^2 - x_2 ^2)/2  $:
	(a) is a trajectory of the original equation \eqref{eq:one} with rapid rotation; (b) is a trajectory of the averaged equation \eqref{eq:averagedequation} including terms up to $O(\varepsilon^3)$; and (c) is a trajectory of the averaged equation without the $O(\varepsilon^3)$ terms.
	The ``wiggles" in (a) are due to the rapid rotation of the potential; the guiding center transformation \eqref{eq:x-transformation} smooths these out while retaining lower-order properties of the orbit, in particular its precession. Comparing (b) and (c), we see that the magnetic term, while only of order $\epsilon^3$, is responsible for this precession.
	All three plots have $ \varepsilon =0.35$, $ x(0) = (1,0) $, $ \dot x (0)=(0,2) $ and $ 0\leq t\leq 430 $. 
	In particular, (c) shows a periodic orbit \emph{over the same time interval as (a) and (b)}, demonstrating the dynamical significance of the magnetic term.
	}
	\label{fig:comparison}
\end{figure}

It is a fundamental fact of nature that a changing electric field creates a magnetic field---for instance,  a rotating electric dipole creates a magnetic field perpendicular to the plane of rotation. Remarkably, an analogous situation arises in mechanics: a changing {\it  mechanical} force field creates a magnetic-like force upon a particle, 
{\it  provided that the change is of high frequency}.  
\begin{figure}[thb]
 	\captionsetup{format=hang}
	\center{  \includegraphics{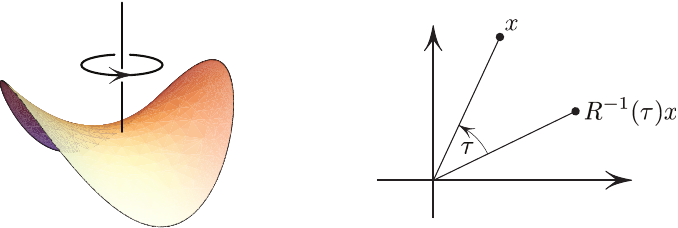}}
	\caption{The value of a   rotating potential evaluated at   a co-rotating point is constant in time: $ U(R( \tau )y, \tau ) = U_0(y) $, 
	and hence $ U(x, \tau ) = U_0(R ^{-1} ( \tau ) x) $.}
	\label{fig:saddlefig}
\end{figure} 
An example,  arising in numerous physical settings, is the motion of a particle in a saddle potential   undergoing rapid rotations, as shown in 
Figure~\ref{fig:saddlefig}. A rotating potential $ U(x,t) $ is obtained by composing a fixed potential $ U_0 (x), \  x\in {\mathbb R}  ^2 $, with the
 rotation through a rapidly changing angle $ t/ \varepsilon $: 
\[
	U: = U_0\circ R^{-1} (t/ \varepsilon ), 
\] 
(see the caption of  Figure~\ref{fig:saddlefig}) and the particle moves according to Newton's second law:
\begin{equation} 
	\ddot x = - \nabla U(x,t/ \varepsilon ),
	\label{eq:one}  
\end{equation}  
where $ \nabla $ is the gradient with respect to $x$. A remarkable phenomenon, known since Brouwer \cite{Brouwer}, is the stabilization of a saddle, such as $U_0(x)= (x_1 ^2 - x_2 ^2)/2 $, under its   rapid rotations (small $\varepsilon$), as illustrated in Figure~\ref{fig:comparison}.  The effective force pulling the particle towards the origin  and responsible for the stabilization is referred to as  the {\it  ponderomotive force} \cite{Schmidt}.  

This restoring force, however, does not account for the precessional motion shown in Figure~\ref{fig:comparison}. It turns out that there is an additional   effective force responsible for this precession.  
This force has the same mathematical expression as the Lorentz magnetic force exerted on a charged particle moving in a magnetic field, with the ``magnetic field" here given in terms of the potential. 
The {\it  ponderomotive Lorentz force} thus seems to be be a fitting name for this magnetic-like force. The existence of this force in the linear case  (more specifically,  for the  quadratic saddle potential rotating with a large constant angular velocity)  was demonstrated in \cite{Kirillov-Levi2}.

The main point of this note is to show that this ponderomotive magnetism  is in fact a ubiquitous phenomenon: it arises generically in potential force fields  with rapid periodic time dependence, including but not limited to rotational dependence. 

More specifically, in the main theorem we assign to any 
solution $ x=x(t) $  of (\ref{eq:one}) with an arbitrary time-periodic potential $ U(x, \tau) $ its ``guiding center" $X=X(t, \varepsilon )$, which extracts the average motion of $x$ by discarding the high-frequency oscillatory part, as shown in Figure~\ref{fig:comparison}. The fundamental observation (and main point) of this paper is that the guiding center $X$ moves as a charged particle in a time-independent potential field 
   {\it  subject to an additional Lorentz magnetic-like force}; see Theorem \ref{thm:ode} for a detailed statement.  The corresponding magnetic vector potential is given explicitly as the time average of  an expression involving   temporal antiderivatives of $U$.  
  
   Although the magnetic effect analyzed in this paper is of one order higher in 
  $\varepsilon$ than the ponderomotive force, it can make a qualitative difference in  the dynamics, as illustrated in  Figures~\ref{fig:comparison} and \ref{fig:precession}.
  
  The ponderomotive force has been well studied, but the ponderomotive {\it   Lorentz} force does not seem to have been explored apart from the linear case of  \cite{Kirillov-Levi2} and a formal treatment of the general case in \cite{S01}.

\subsection*{Some related background}
Almost a century has passed since Brouwer analyzed the motion of a particle in a rotating saddle potential  \cite{Brouwer}. Brouwer considered  a particle sliding without friction on a rotating saddle surface in the presence of gravity,   and analyzed the linearized equations near the equilibrium. One consequence of this analysis was the conclusion that the equilibrium is linearly stable if the two principal curvatures are equal and the rotation is sufficiently rapid. The rotating saddle problem arises in numerous settings in physics, a long list of which can be found in \cite{Kirillov2003} and \cite{Kirillov-Levi1}. We mention only one example here: the motion of asteroids near the equilateral Lagrangian points $ L4$ and $ L5$. The same old problem acquired new relevance starting in the 1990s due to  physical applications---a rotating saddle potential can be used to trap charged particles, for instance---and it has since been studied more extensively, both analytically and numerically \cite{HB2005,Kirillov2011}.

Stabilization caused by rapid rotation of the saddle is, 
at least at the first glance, quite counterintuitive, although it admits a simple heuristic explanation.

A related   phenomenon which shares this counterintuitive defiance of gravity with the rotating saddle is the stabilization of an inverted pendulum through rapid vibration of its pivot; this was discovered more than a century ago (1908) by Stephenson \cite{stephenson} but is often referred to as the Kapitsa effect, having been studied by Kapitsa \cite{Kapitsa1951} in the 1950s. In fact, the inverted pendulum was a motivation for the invention of the Paul trap used to confine charged particles, for which W. Paul was awarded the 1989 Nobel Prize in physics \cite{Paul1990}. The Kapitsa effect also has unexpected connections to differential geometry \cite{Levi1998} and non-holonomic systems  \cite{Weckesser}.

In this paper we extend the above mentioned results by proving a general averaging theorem which accounts for the effects of rapid periodic changes in a potential force field acting on a particle. According to this theorem, the particle is subject to a magnetic-like force in addition to an effective potential, the latter resulting in the ponderomotive force described above. 

\subsection*{Physical applications and motivation} A surprisingly large class of physical problems is described by   (\ref{eq:one}). 
These include the dynamics of electron beams in helical quadrupole magnetic
fields (the motion in the direction of the helical axis results in the rotation of the potential in the 
transversal plane, so that the  displacement along the axis plays the role of the time);  the dynamics of ions in rotating ion traps;  the restricted circular three-body problem in celestial mechanics,
 describing the motion near the triangular Lagrange libration points L4 and L5; and much more. A more extensive (but by no means comprehensive) list can be found in \cite{Kirillov-Levi1}.
Ponderomotive magnetism arises in all of these problems.  
 
\section{Ponderomotive magnetism by averaging}
\subsection*{The setting and the notations}

We consider the motion of a particle in a rapidly oscillating potential,
 \begin{equation} 
	\ddot x = - \nabla U (x, t/ \varepsilon ) , \  \  x \in {\mathbb R}  ^n , 
	\label{eq:mainequation}
\end{equation} 
where the gradient is taken with respect to $x$, the function $U$ is periodic in the second variable:
\[
	U(x, \tau)=U(x, \tau+1), 
\]   
and $\varepsilon$ is a small parameter. This general form covers the two special cases of rotating and oscillating potentials mentioned in the introduction. Note that \eqref{eq:mainequation} can be written as a Hamiltonian system, with the Hamiltonian
\begin{equation} 
	H(x,p,t/\eps) = \frac12 p^2 + U(x,t/\eps).
	\label{eq:hamiltonian}
\end{equation}   
 
 To state our result, we introduce the following temporal antiderivatives of $ U $: 
\begin{equation} 
	V(x, \tau ) = \int [U(x, \tau )-\overline {U}(x)] \,d \tau, \  \  S(x, \tau ) = \int V(x, \tau ) \,d \tau, \  \ A(x, \tau ) = \int S(x, \tau ) \,d \tau, \  \ 
	\label{eq:VSA}
\end{equation}   
where $\overline U = \overline {U}(x) $ is the time average of $U$ over one temporal period, and where each of the indefinite integrals is chosen to have zero time average:
\[
	\overline V = \overline S = \overline A = 0. 	
\] 

\begin{rem}    If $S$ is interpreted as the position of a particle at time $t$, then the antiderivative $A$ is called the \emph{absement}. Physically, the absement can be thought of as follows: If $S$ is the distance by which the gas pedal of a car is depressed, then the absement is proportional to the total amount of fuel provided, and hence the total amount of energy available  to the car.\footnote{We thank Steve Mann for helpful comments regarding physical interpretations of the absement.}
\end{rem}

\noindent {\bf Notation.} 
\label{notation} From now on, for the sake of brevity, the $x$-derivative of a function $ f\colon {\mathbb R}^n \rightarrow {\mathbb R}$ will be denoted by a prime $\ ^\prime$ rather than by $ \nabla $. Similarly,  the second derivative, i.e. the Hessian,  will be denoted by 
two primes, so that  
\begin{equation*} 
	U^\prime = \nabla U, \  \  U^{\prime\prime} = \nabla ^2  U,  \  \  \hbox{etc.} 
\end{equation*}
 
 \subsection*{The main results}
 
The main results of this paper are the following two theorems, the first of which describes the averaging in the phase space.

\begin{thm} \label{thm:1} Assume that the potential $U$  in \eqref{eq:mainequation} is of class $C^3$ in $x$ and continuous in $t$. For any fixed ball in the {\rm phase space} ${\mathbb R} ^{2n}$
there exists $\varepsilon_0>0 $  such that for each positive  $\varepsilon < \varepsilon_0$ there is a time-parametrized family of symplectic transformations 
\begin{equation*}
	\psi^\tau =\psi^\tau_ \varepsilon\colon (x,p) \mapsto (X,P)    
\end{equation*}
defined on the ball and periodic of period one in $\tau$, such that the time-dependent transformation 
  $ \psi^{t/ \varepsilon}  $ 
 turns the system with the Hamiltonian \eqref{eq:hamiltonian} into 
the system with the Hamiltonian 
\begin{equation} 
	K(X, P, t/\eps) = \frac12 P^2 + \overline{U} + \frac{\eps^2}{2} \overline{V' \cdot V'} - \eps^3 \big(\overline{S'' V'} \cdot P\big) + O(\eps^4), 
\label{eq:avhamiltonian}
\end{equation}   
which is time independent to third order.
\end{thm} 

The second theorem  describes the averaging in the configuration space.  For the Newtonian formulation \eqref{eq:mainequation} we have the following, recalling the notation from \eqref{eq:VSA}.

\begin{thm} \label{thm:ode} Assume that $U$ is of class $C^4$ in $x$ and $ C^1 $  in $t$. 
Let   $x=x(t)$ be a solution  of \eqref{eq:mainequation} which stays in a fixed ball in $ {\mathbb R}  ^n $ for a time interval\footnote{For many potentials, e.g. superquadratic ones such as $ x^4-y^4 $,  solutions may blow up in finite time.} $I$. There is an associated function of the form
\begin{equation} 
	X =  x + \eps^2 S'(x,t/ \varepsilon ) - 2\eps^3 A''(x,t/ \varepsilon) \dot x + O(\varepsilon ^4)
	\label{eq:x-transformation}
\end{equation}  
that satisfies the equation
\begin{equation} 
	\ddot X = - \overline {U}\, ^\prime(X) - \varepsilon ^2 W ^\prime(X) + \varepsilon ^3 B(X) \dot X + O ( \varepsilon ^4),  
	\label{eq:averagedequation}
\end{equation}   
for all $ t\in I $, where $W=\frac{1}{2} \overline{V ^\prime \cdot V ^\prime}$ and $B$ is the skew-symmetric matrix given
by $B = (b ^\prime) ^T-b ^\prime$, with $b = \overline {S ^{\prime\prime} V ^\prime} $.  
\end{thm}

We refer to $X$ as the \emph{guiding center}\footnote{although the term is commonly used in a different context: that of motion of a charged particle in magnetic field.} associated to $x$, and to \eqref{eq:averagedequation} as the \emph{averaged equation}. We do not specify the $O(\varepsilon^4)$ term in \eqref{eq:x-transformation} for the sake of simplicity; it can easily be obtained from the proof if desired.

Theorem \ref{thm:ode} is an extension of the main result in \cite{Kirillov-Levi2} to the nonlinear case and allows for arbitrary periodic time dependence (in particular, the potential need not be rotational).     An attempt at a straightforward extension of the method used there quickly becomes too cumbersome, and is probably close to impossible without a computer. We achieve a drastic simplification of the proof by using the Hamiltonian nature of the system. The proof is given in Section \ref{sec:proofs1}.

  \begin{figure}[thb]
 	\captionsetup{format=hang}
	\center{  \includegraphics{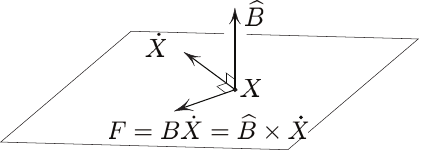}}
	\caption{$ b = \overline{S ^{\prime\prime} V^\prime }$ is the vector magnetic potential: in dimension $ n=3 $, 
	$   \widehat B \buildrel{def}\over{=} {\rm curl}\,  b $ is the magnetic field. In dimension $ n = 2 $, the scalar $ {\rm curl}\, b $ is the magnitude of the  magnetic field normal to the plane.  }
	\label{fig:magforce}
\end{figure}  
Before giving the proof of this theorem, we make several comments. 
\begin{enumerate} 
\item   {\bf The effective Lorentz force.} The term $ \varepsilon ^3 B(X) \dot X $ in  \eqref{eq:averagedequation} can be interpreted as the force produced by a magnetic field upon a charge $ \varepsilon ^3 $, as in Figure~\ref{fig:magforce}. In dimension $ n = 3 $ one has $ B(X) \dot X = \widehat B \times \dot X  $, where 
the magnetic vector   $ \widehat B= (\widehat B_1, \widehat B_2, \widehat B_3) $ is associated to the matrix $B$ via the standard isomorphism between $ {\mathbb R}  ^3$ and  the Lie algebra $ so(3) $: 
\[
	B= \left( \begin{array}{ccc} 0 & -\widehat B_3 & \widehat B_2\\   \widehat B_3 & 0 & -\widehat B_1    \\ -\widehat B_2 & 
	\widehat B_1 & 0\end{array} \right) .
\]  
Thus $ \varepsilon ^3 B(X) \dot X  $ looks exactly like the standard expression for the Lorentz magnetic force, $F = q {\bf v} \times {\bf B}$,
with $ q = \varepsilon ^3 $, $ {\bf v} = \dot X $ and ${\bf B} = - \widehat B $.

\item {\bf The magnetic vector potential.} In dimension $n = 3 $ the magnetic field is $ \widehat B = {\rm curl}\,   b $. In other words,   $b = \overline{S ^{\prime\prime} V^\prime }$ 
 is the {\it  magnetic vector potential of the magnetic field. }  
In dimension $ n =2 $, the strength of the magnetic field normal to the plane is given by $ {\rm curl}\, b $,  the scalar curl in $ {\mathbb R}  ^2 $. 

 \item {\bf The ponderomotive Lorentz force causes precession.} Consider a rotating potential in dimension
 $ n=2 $. Since the averaged potential $W$ in this case is rotationally symmetric,  every  trajectory of the averaged system {\it  without the magnetic term} that passes through the origin  simply oscillates along a straight line segment, as in Figure~\ref{fig:precession}(a). 
Although the magnetic term is not  of leading order in the averaged equation, it is responsible for the precession as illustrated in Figure~\ref{fig:precession}(b) and (c); see also Figure \ref{fig:comparison}. The same dominant effect of precession will manifest itself for trajectories passing near the origin. It should be noted that  trajectories not passing through the origin for generic  (rotationally symmetric) potentials do precess without the magnetic effect; this precession may dominate the effect of ponderomotive magnetism. However, the ponderomotive magnetism does have a major effect (in the sense of Figure~\ref{fig:precession}) for those  trajectories that pass through or near the origin, i.e. for the ones with small angular momenta relative to their energy.

 \begin{figure}[thb]
 	\captionsetup{format=hang}
	\center{  \includegraphics{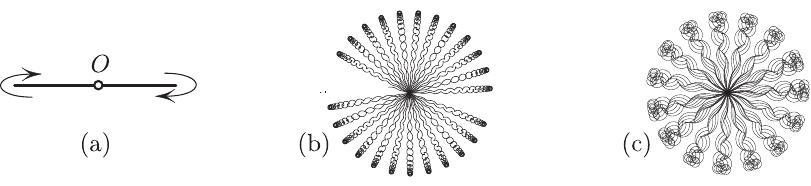}}
	\caption{Precessional effect of the ponderomotive Lorentz force: (a) is a trajectory of the averaged system   {\it  without} the magnetic term, while (b) and  (c) are  trajectories through the origin of the averaged equation with magnetic term included, for $ \varepsilon = 0.35 $ and $\varepsilon = 0.5 $, respectively. Here the potential is $ U_0(x,y) = \cos x- \cos y $, the initial conditions are $ x(0)=y(0)= \dot y (0)=0 $, $ \dot x (0)=0.2 $, and the time interval is $ 0\leq t \leq 700 $.}
	\label{fig:precession}
\end{figure}

 \item {\bf Linearity in the velocity.} Although the original system \eqref{eq:mainequation} is fully nonlinear, its averaged counterpart  \eqref{eq:averagedequation} is, to leading order, linear in the velocity $\dot X$.

 \item {\bf Possible blow-up.} Unlike in the case of quadratic potentials, which correspond to linear systems, solutions may in general blow up in finite time. An interesting open problem is 
to understand the nature of the set of initial data for which solutions avoid blow-up.  In the case of rotating potentials the problem reduces to studying time-independent potentials with added Coriolis force. We do not discuss this problem since it is peripheral to the main point of this 
note. 
\item {\bf An adiabatic invariant.} As a side remark we note that the rapid  time dependence in the original system   (\ref{eq:mainequation}) gives rise to an approximately conserved quantity. Indeed, for the averaged system   (\ref{eq:averagedequation})  the energy   
\begin{equation} 
	E_{\rm av}(X, \dot X) = \frac{\dot X ^2 }{2} +\overline U(X) + \varepsilon^2 W(X) 
	\label{eq:energy}
\end{equation}   
 is nearly conserved in the sense that along solutions confined to a bounded region one has
\[
	\frac{d}{dt} E_{\rm av}(X, \dot X) = O ( \varepsilon ^4),
\]  
as is easily verified. Now to any solution  $x$ of   (\ref{eq:mainequation}) we assign the ``effective energy" defined by 
  $ E(x, \dot x ) \buildrel{def}\over{=}  E_{\rm av}(X, \dot X)$, where $x$ and $X$ are related via   (\ref{eq:x-transformation}), so that
\[
	\frac{d}{dt} E(x, \dot x ) = O( \varepsilon ^4 )
\]   
along solutions of   (\ref{eq:mainequation}) (which are likewise confined to a bounded region).  For small $\varepsilon$ this implies  small change in $E$  over a long time interval. With some natural assumptions this implies $ O( \varepsilon ^ {4- \alpha} )$ change of $E(x, \dot x ) $ over the time interval $ 0\leq t\leq  \varepsilon^ {-\alpha} $ for $ 0< \alpha < 4 $. We leave out all the details since this question lies outside the focus of this paper.

\item {\bf An additional adiabatic invariant   in the rotational case and near-integrability.}  For a rotating potential in dimension $ n=2 $, namely  $ U(x,t) =U_0 ( R(-t/ \varepsilon )x )$, both averaged  potentials $ \overline U $ and $ W $, as well as   the magnetic vector 
potential $ b = \overline{S ^{\prime\prime} V ^\prime } $, are rotationally symmetric, i.e. they depend only on $ | X | $.  In this case the  truncated system has a conserved quantity associated with this symmetry, namely  
\begin{equation} 
	I_{\rm av}( X, \dot X) \buildrel{def}\over{=} X\wedge \dot X +   \varepsilon ^3 \int_{0 }^{| X |} r \beta (r) dr, 
	\label{eq:am}
\end{equation}  
where $ \beta = {\rm curl}\,  b$. The first term in this sum is the angular momentum of a unit point mass at $X$;  the integral in   (\ref{eq:am})   is the 
flux of the magnetic field through the disk of radius $ | X | $ (up to the factor of $ 1/2 \pi $). As before,     to any solution $ x$ of the original system   (\ref{eq:mainequation}) we associate the function  
\[
	I (x, \dot x ) \buildrel{def}\over{=} I_{\rm av}( X, \dot X),
\]  
where $X$ is related to $x$  via   (\ref{eq:x-transformation}).  With this definition, 
$ \frac{d}{dt} I(X, \dot X) = O ( \varepsilon ^4) $, and the remark at the end of the preceding paragraph applies.  Furthermore,  the functions 
 $ E_{\rm av} $ and $ I_{\rm av} $  expressed in terms of position and momentum are in involution. Thus the Hamiltonian form of  (\ref{eq:averagedequation})  is a small (in a certain sense) perturbation of a completely integrable Hamiltonian system  corresponding to 
   \eqref{eq:avhamiltonian}  with the $ O(\varepsilon ^4)$ terms deleted. It would be of interest to see whether one can apply KAM theory to show that the  extended phase space of   (\ref{eq:averagedequation}),   and thus   of    (\ref{eq:mainequation}), has   invariant tori whose combined measure approaches full measure as $ \varepsilon \rightarrow 0 $ (with appropriate assumptions guaranteeing in particular the avoidance of blow-up). This would imply that most (in the Lebesgue sense) initial data lead to solutions bounded for all time, although one would still expect Arnold diffusion for a small (in the Lebesgue sense) set of solutions.

\end{enumerate}

In the next section we discuss some implications and applications of Theorem \ref{thm:ode}.

\section{Applications of the averaging theorem}\label{sec:examples}
The examples of this section give particularly simple illustrations of   Theorem \ref{thm:ode}. These are:
\begin{enumerate} 
\item Rotating separable potentials;

\item Purely oscillatory potentials; 
\end{enumerate}    

These two classes encompass a broad range of physical systems---such as the rotating saddle, the Paul trap, and the inverted pendulum with oscillating pivot---and yet are simple enough that the magnetic term in the averaged equation can be worked out explicitly.

\subsection*{The rotating separable potential}
By a rotating separable potential we mean a rotating planar potential $U(x,\tau) = U_0(R ^{-1} (\tau) x)$, as described in the introduction, with the additional requirement that $U_0$ can be written in polar coordinates as $U_0(x) = f(r) h(\theta)$, where $x = (r \cos\theta, r \sin\theta)$. In this case the time-dependent potential becomes 
\begin{equation}\label{Usep}
	U(x,\tau) = f(r) h(\theta - \tau),
\end{equation}
where $h$ is a $2\pi$-periodic function. Note that for any function of $r$ and $\theta - \tau$, say $k(r,\theta - \tau)$, that is $2\pi$-periodic in the second variable, the time average satisfies
\begin{equation}\label{eq:kavg}
	\bar k(r,\theta) = \frac{1}{2\pi} \int_0^{2\pi} k(r,\theta - \tau) \, d\tau
	= \frac{1}{2\pi} \int_0^{2\pi} k(r,\tau) \, d\tau
\end{equation}
for each $\theta$. That is, \emph{the time average $\bar k$ depends only on $r$.} This fact greatly simplifies the following computation, since it tells us that certain terms in the averaged equation will vanish, and hence do not need to be calculated explicitly.

We now work out the averaged equation, as given by Theorem \ref{thm:ode}. The required temporal antiderivatives are
\begin{align*}
	V(x,\tau) = -f(r) v(\theta - \tau), \quad S(x,\tau) = f(r) s(\theta - \tau),
\end{align*}
where $v$ and $s$ are $2\pi$-periodic functions satisfying $v' = h$, $s' = v$, and $\bar v = \bar s = 0$. Letting $\ber = (\cos\theta, \sin\theta)$ and $\bet = (-\sin\theta, \cos\theta)$ denote the standard polar unit vectors, we have
\[
	V' = -f' v \,\ber - \frac{f h}{r} \bet
\]
(with all terms evaluated at $r$ and $\theta - \tau$) and hence
\begin{equation}\label{eq:Wsep}
	W = \frac12 \overline{V' \cdot V'} = \frac12 \left[ (f')^2 \overline{v^2} + \frac{f^2}{r^2} \overline{h^2} \right].
\end{equation}

We next compute the angular component of the vector field $b = \overline{S''V'}$. (According to \eqref{eq:kavg}, the radial component of $b$ will be of the form $g(r) \ber$ for some function $g$; therefore its curl vanishes and it does not contribute to the magnetic field $B$ in \eqref{eq:averagedequation}.) Using

\[
	S' = f' s \,\ber + \frac{f v}{r} \bet
\]
and the identities $\frac{\p}{\p\theta} \ber = \bet$,  $\frac{\p}{\p\theta} \bet = -\ber$, we compute
\[
	S'' \ber 
	= \frac{\p}{\p r} S' = \frac{r f' - f}{r^2} v \bet + (\ast) \ber
\]
where $ (\ast)  $ denotes a function of $r$ and $\theta - \tau$, and 
\[
	S'' \bet = \frac{1}{r} \frac{\p}{\p\theta} S' = \frac{1}{r} \left( f' s + \frac{f h}{r} \right) \bet + (\ast) \ber,
\]
from which we find
\begin{align*}
	S''V' = -f'v (S'' \ber) - \frac{fh}{r} (S'' \bet) = \left[ -f'v \left(\frac{rf' - f}{r^2} \right) v - \frac{fh}{r^2} \left( f's + \frac{fh}{r} \right) \right] \bet + (\ast) \ber.
\end{align*}
Finally, using the relation $\overline{hs} = \overline{v's} = - \overline{s'v} = - \overline{v^2}$, we obtain the magnetic vector potential
 \begin{equation}\label{eq:bsep}
 	b = \overline{S''V'} = \left[ \frac{2ff' - r (f')^2}{r^2} \overline{v^2} - \frac{f^2}{r^3} \overline{h^2} \right] \bet + (\ast) \ber,
\end{equation}
where now the coefficient of $\ber$ only depends on $r$, by \eqref{eq:kavg}.

We now examine some special cases using \eqref{eq:Wsep} and \eqref{eq:bsep}. Throughout we assume that $\bar h = 0$, so that $\overline U = 0$ and hence there is no zeroth-order term in the averaged equation \eqref{eq:averagedequation}.

\begin{rem}
For the magnetic term to dominate \eqref{eq:averagedequation} we must have $W$ constant (so that $W'$ vanishes). From \eqref{eq:Wsep} we see that this is only possible when $f(r) = Ar$ for some constant $A$. Then \eqref{eq:bsep} implies that the angular component of $b$ is proportional to $r^{-1} \bet$, and so $\operatorname{curl} b = 0$. We conclude that there are no potentials of the form \eqref{Usep} for which the magnetic term is leading order.
\end{rem}

\begin{figure}[thb]
 	\captionsetup{format=hang}
	\center{  \includegraphics{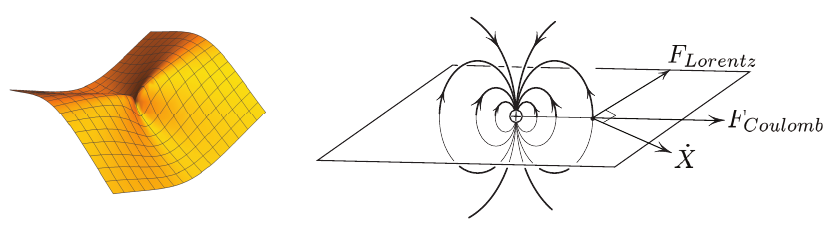}}
	\caption{Rotation of the potential (\ref{eq:example}) (the graph is shown on the left) gives rise to a ``repelling charge" and a  ``magnetic dipole" acting on the guiding center with forces  $F_{Coulomb}$ and $F_{Lorentz}$ lying in the plane of motion. }
	\label{fig:ruledsurface}
\end{figure}

Suppose $f(r) = r^m$ for some real number $m$, so the potential is homogeneous. In this case we find
\[
	W = \frac12 \left( m^2 \overline{v^2} + \overline{h^2} \right) r^{2m-2}
\]
and
\[
	b = \left( (2m - m^2) \overline{v^2} - \overline{h^2} \right) r^{2m-3} \bet + (\ast) \ber
\]
In particular, the magnitude of $W'$ is proportional to $r^{2m-3}$ and the magnitude of $B$ is proportional to $r^{2m-4}$. Recalling the form of the averaged equation
\[
	\ddot X = - \overline {U}\, ^\prime(X) - \varepsilon ^2 W ^\prime(X) + \varepsilon ^3 B(X) \dot X + O ( \varepsilon ^4),  
\]
we conclude that \emph{the magnetic term will dominate $\varepsilon^2 W'(X)$ when $X$ is small relative to $\varepsilon$.} 
In particular, choosing $m = 1/2$ we find that $W$ and $b$ are proportional to $1/r$ and $1/r^2$, respectively. That is, $W$ is equal to  the potential of a point change, and $b$ is equal to the vector potential of a magnetic monopole with axis normal to the plane; see Figure~\ref{fig:ruledsurface}. For such a potential, assuming that $ \bar h = 0 $, the averaged equation is 
\begin{equation}
	\ddot X =\underbrace{  \varepsilon^2 \alpha\frac{   X}{  |X|^3}}_{F_{\rm Coulomb}} \ -\  \underbrace{2 \varepsilon ^3 \beta \frac{J \dot X}{  |X|^3}}_{F_{\rm Lorentz}}  + O(\varepsilon^4), \qquad J = \begin{pmatrix} 0 & -1 \\ 1 & 0 \end{pmatrix}, 
\end{equation}
with $\alpha$ and $\beta$ being the coeffients in $W$ and $b$. 
An example of such a potential, shown in Figure~\ref{fig:ruledsurface}, is
\begin{equation}
	U(x) = \sqrt{r} \cos 2\theta = \frac{x^2 - y ^2 }{(x^2 + y^2)^{3/4}},
	\label{eq:example}
\end{equation}
which has  $ \alpha = 17/64$ and $\beta = -13/32 $.

Another case of great physical importance, as described in \cite{Kirillov2003}, is the rotating quadratic saddle potential, i.e. $U_0(x ) = (x_1 ^2 - x_2 ^2 )/2$.
This is of the general form \eqref{Usep}, with $f(r) =r^2$ and $h(\theta) = \frac12 \cos 2\theta$. It follows that
\[
	\bar h = 0, \quad \overline{h^2} = \frac{1}{8}, \quad \overline{v^2} = \frac1{32}
\]
and so
\[
	W = \frac{r^2}{8}, \qquad b = - \frac{r}{8} \bet + (\ast) \ber,
\]
from which we obtain the averaged equation
\begin{equation}
	\ddot X = - \frac{\varepsilon^2}{4} X + \frac{\varepsilon^3}{4} J \dot X + O(\varepsilon^4),
	\label{eq:avlin}
\end{equation}
in agreement with the result in \cite[Theorem 1]{Kirillov-Levi2}.

 \subsection*{The purely oscillatory potential}
\begin{figure}[thb]
 	\captionsetup{format=hang}
	\center{  \includegraphics{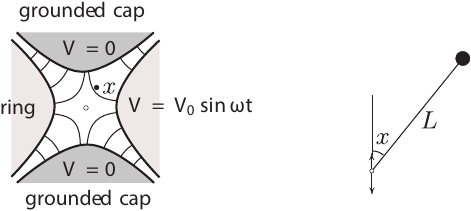}}
	\caption{Left: Axial cross-section of a Paul trap: a ring with two caps, one above and one below. 
	 Right: inverted pendulum with a vertically vibrating pivot: a point mass on a weightless rod.  }
	\label{fig:paultrap}
\end{figure}

 Another special case of physical importance is 
 \begin{equation} 
	\ddot x  = -a(t/ \varepsilon ) U ^\prime (x),  	
	\label{eq:oscillatory}
\end{equation} 
where $a$ is real-valued and periodic. This equation describes, among other things, the motion of charged particles in the Paul trap \cite{Paul1990}, in which an electric field is generated by oscillating voltages on electrodes, as shown in Figure~\ref{fig:paultrap} (left). The 
position $x$ of a 
charged particle in the cavity is then governed by    (\ref{eq:oscillatory}). (If gravity is included, an additional constant term should be added to the right-hand side of   (\ref{eq:oscillatory}).) With an oscillating  voltage  applied to the ring, the electostatic field in the cavity looks like the right-hand side of  (\ref{eq:oscillatory}). 
	This oscillation can stabilize the equilibrium.  By contrast, in a static electric field this equilibrium is necessarily 
	unstable since harmonic potentials cannot have minima.

The  inverted pendulum with oscillating suspension    (Figure~\ref{fig:paultrap}, right)  is also governed by   (\ref{eq:oscillatory}). In this case $ a(t) =  -L ^{-1}(g+ a_{\rm vert} (t) )$, where $ a_{\rm vert} $  is the vertical acceleration of the pivot and $U(x) =- \cos x $.  

For \eqref{eq:oscillatory}, the averaged equation \eqref{eq:averagedequation} becomes 
\[
	\ddot X=-\overline {a}\,U ^\prime (x) -  \varepsilon ^2\, \overline{v ^2} \,U ^{\prime\prime}(x) U ^\prime (x) + O( \varepsilon ^4), 
\] 
where $v(\tau) = \int [a(\tau) - \bar a] \,d\tau$ is chosen to have $ \overline v = 0 $. Note in particular that the $O(\varepsilon^3)$ terms vanish and so there is no magnetic term in this case. This adds new insight into the behavior of   (\ref{eq:oscillatory}), showing that, at least to third order in $\varepsilon$, there is no precession. This question was not addressed in earlier studies, e.g. \cite{Levi1999}, which did not compute cubic terms.  This example also suggests that the $ O ( \varepsilon ^3 ) $ magnetic term is caused by rapid rotation of the potential, rather than oscillation.

\section{Proof of Theorems \ref{thm:1} and \ref{thm:ode} }

\subsection{The general pattern}
Before proceeding with the proof, we describe the recurring inductive step in the normal form reduction. This step starts with  the Hamiltonian
\begin{equation} 
	{\cal H}_k (x, p, t/ \varepsilon ) =H_0+ \varepsilon H_1 +\cdots + \varepsilon ^{k-1}H_{k-1}+  
	\varepsilon ^{k} \tilde H_{k}+\varepsilon ^{k+1} \tilde H_{k+1}+\cdots, 
	\label{eq:H}
\end{equation}   
where the time dependence has been eliminated in the terms without the tilde: $   H_j =   H_j(x, p) $ for $ j\leq k-1 $, while     
$  \tilde H_{k } =  \tilde H_k(x,p, t/ \varepsilon ) $. 
Wishing to eliminate the time dependence in $ \tilde H_k $ by a time-periodic symplectic transformation $(x,p)\mapsto(x ^\prime ,p^\prime ) $, we  seek the generating function $ G (x,p, \tau ) $, periodic in $\tau$ of period $1$,  for such a transformation: 
\begin{align}
	\begin{split}
	x^\prime  &= x \ + \eps^{k+1} G_p(x, p^\prime ,t/\eps) \\
	p &= p^\prime  + \eps^{k+1} G_x(x,p^\prime ,t/\eps).   
	\end{split}\label{eq:gs}
\end{align}
If $G$ is bounded in the $C^1$ norm in $ (x, p ^\prime) $, then by the implicit function theorem \eqref{eq:gs} defines a map $ \varphi \colon (x,p)\mapsto(x ^\prime ,p^\prime )$, depending on the parameter $t$, and $ \varphi ^{-1} $   is  given by   
\begin{align}
	\begin{split}
	x  &= x^\prime - \eps^{k+1} G_p(x ^\prime ,p ^\prime ,t/\eps) + O(\varepsilon ^{2k+2})\\
	p   &= p ^\prime   + \eps^{k+1} G_x(x ^\prime ,p ^\prime ,t/\eps)+O(\varepsilon ^{2k+2}),  
	\end{split}\label{eq:gs1}
\end{align}
where the remainder is a function of $x^\prime$, $p ^\prime$ and $ \tau = t/ \varepsilon  $. The $ C^1 $-boundedness of $G$ will be established in the course of its construction. The transformed Hamiltonian is given by 
\begin{equation} 
	 {\cal H}_{k+1} ( x^\prime , p ^\prime , t/ \varepsilon ) =   {\cal H}_k( x, p, t/ \varepsilon  ) + 
	\varepsilon  ^{k+1} \partial _t G(x, p ^\prime , t/ \varepsilon )  , 
	\label{eq:newham}
\end{equation} 
where  $x$ and $p$ are functions of the new independent variables $ x ^\prime $ and $ p^\prime $ by virtue of   \eqref{eq:gs1}, and where 
\[
	 \partial_t G(x, p ^\prime , t/ \varepsilon ) = \varepsilon ^{-1} \partial_\tau G(x,p,\tau) \  \  \hbox{with}  \  \   \tau =  t/ \varepsilon , 
\] 
so that  the correction term in  (\ref{eq:newham})  is of order $ \varepsilon ^k $: 
\begin{equation} 
	 {\cal H}_{k+1} ( x^\prime , p ^\prime , t/ \varepsilon ) =   {\cal H}_k( x, p, t/ \varepsilon  ) + 
	\varepsilon  ^{k} \partial _\tau G(x, p ^\prime , t/ \varepsilon ). 
	\label{eq:newham1}
\end{equation} 
Substituting \eqref{eq:H} and \eqref{eq:gs1} into \eqref{eq:newham1}, we obtain the transformed Hamiltonian
\begin{equation} 
	 {\cal H}_{k+1}   = H_0  +  \cdots + 
	 \varepsilon ^{k-1}H_{k-1} + \varepsilon ^k \big(\tilde H_k +
	 \partial_\tau G  \big) +O ( \varepsilon ^{k+1} ),  
	\label{eq:newham2}
\end{equation}   
where all the functions are now evaluated at the new variables $ ( x ^\prime , p ^\prime ) $. To complete the inductive step, we make $ \tilde H_k + \partial_\tau G $ time independent by setting $ G = - \int (\tilde H_k-\overline {\tilde H}_k)\,d\tau $, with the constant of integration chosen so as to have $ \overline G = 0$.

Note that the terms of order $ \leq k-1 $ are not affected by the transformation. The Poisson bracket $ \{H_0, G\} $ appears  in the $ \varepsilon ^{k+1} $-term, while $ \partial_\tau G$ appears in the $ \varepsilon ^k$-term; this mismatch is due to the rapid time dependence of the Hamiltonian. When implementing this general scheme we will have to keep track of the higher-order terms, as they will be averaged in subsequent steps.

\label{sec:proofs1}
\subsection{Proof of Theorem  \ref{thm:1}}
We start with the Hamiltonian
\begin{align}
	H(x,p,t/\eps) = \frac12 p^2 + U(x,t/\eps).
\end{align}
To eliminate time dependence in the $ \varepsilon ^0 $-term $U$, consider a generating function $G^{(1)}(x, p, \tau )$ of period $1$ in $\tau$ (the superscript refers to the first step in the averaging). Such a function   defines a symplectic map  $(x,p) \mapsto (x_1,p_1)$  implicitly via
\begin{align*}
	x_1 &= x + \eps G^{(1)}_p(x,p_1,t/\eps) \\
	p &= p_1 + \eps G^{(1)}_x(x,p_1,t/\eps),  
\end{align*}
and the transformed system is Hamiltonian with the new Hamiltonian $ H_1=H+\frac{\partial}{\partial t}  (\varepsilon G^{(1)})=H+G^{(1)}_\tau $, 
where $ G^{(1)}_\tau $ denotes the derivative with respect to the last variable: 
\[
	H_1(x_1,p_1,t/\eps) = H(x,p,t/\eps) +   G^{(1)}_\tau (x, p_1, t/ \varepsilon). 
\]
Recalling the form of $H$, we have 
\[
	H_1 = \frac12 p_1^2 + U(x,t/\eps) + G^{(1)}_\tau + O(\eps).
\]
We now choose $G^{(1)}$ so as to cancel the time dependence in $U$: 
\[
	-G^{(1)}(x,p_1,\tau) = V(x,\tau) := \int [U(x,\tau) - \overline{U}(x)]d\tau,
\]
with the constant of integration chosen so that $\overline{V} = 0$. Since such $G^{(1)}$ is independent of $p_1$, the transformation simplifies to
\begin{align}\label{cov1}
\begin{split}
	x_1 &= x \\
	p &= p_1 - \eps V'(x,t/\eps),
\end{split}
\end{align}
and so the new Hamiltonian is
\[
	H_1 = \frac12 p_1^2 + \overline{U}(x_1) - \eps p_1 \cdot V'(x_1,t/\eps) + \frac{\eps^2}{2} (V' \cdot V')(x_1,t/\eps).
\]

To average the $O(\eps)$ terms, i.e. make them time independent,  we repeat the above transformation, setting
\begin{align*}
	x_2 &= x_1 + \eps^2 G^{(2)}p(x_1,p_2,t/\eps) \\
	p_1 &= p_2 + \eps^2 G^{(2)}_x(x_1,p_2,t/\eps)
\end{align*}
for some new generating function $G^{(2)}$ to be determined.  Since the new Hamiltonian will have the form
\[
	H_2 = H_1 + \eps G^{(2)}_\tau = \frac12 p_2^2 + \overline{U}(x_2) - \eps p_2 \cdot V'(x_1,t/\eps) + \eps G^{(2)}_\tau + O(\eps^2),
\]
we choose $G^{(2)}(x_1,p_2,t/\eps) = p_2 \cdot S'(x_1,t/\eps)$, where $S(x_1,\tau) = \int V(x_1,\tau) d\tau$ is chosen to have zero mean. Therefore the transformation is
\begin{align}
\begin{split}
	x_2 &= x_1 + \eps^2 S'(x_1,t/\eps) \\
	p_1 &= p_2 + \eps^2 S''(x_1,t/\eps) p_2.
\end{split}
\end{align}
To find $H_2$, we first compute
\[
	x_1 = x_2 - \eps^2 S'(x_2,t/\eps) + O(\eps^4)
\]
and
\[
	p_1^2 = p_2^2 + 2 \eps^2 p_2 \cdot \big(S''(x_2,t/\eps) p_2 \big) + O(\eps^4).
\]
Therefore
\begin{align*}
	H_2 &= H_1 + \eps G^{(2)}_\tau \\
	&= \left[ \frac12 p_2^2 + \eps^2 p_2 \cdot \big( S''(x_2,t/\eps) p_2 \big) \right]
	+ \left[ \overline{U}(x_2) - \eps^2 \overline{U}'(x_2) \cdot S'(x_2,t/\eps) \right]  \\
	&\quad\ - \eps p_1 \cdot V'(x_1,t/\eps) + \eps p_2 \cdot V'(x_1,t/\eps) 
	+ \frac{\eps^2}{2} (V' \cdot V')(x_1,t/\eps) + O(\eps^4),
\end{align*}
which simplifies to 
\[
	H_2 = \frac12 p_2^2 + \overline{U}(x_2) + \eps^2 \left[ p_2 \cdot (S'' p_2) - \overline{U}' \cdot S' + \frac12 V' \cdot V' \right] - \eps^3 (S'' V') \cdot p_2
	+ O(\eps^4),
\]
with all terms on the right-hand side evaluated at $x_2$ and $t/\eps$. Note that the $O(\eps)$ terms have averaged to zero.

For the next averaging step we set
\begin{align}
	\begin{split}
	x_3 &= x_2 + \eps^3 G^{(3)}_p(x_2,p_3,t/\eps) \\
	p_2 &= p_3 + \eps^3 G^{(3)}_x(x_2,p_3,t/\eps).
	\label{eq:3}
	\end{split}
\end{align}
Since
\[
	H_3 = H_2 + \eps^2 G^{(3)}_\tau= \frac12 p_3^2 + \overline{U}(x_3) + \eps^2 \left[ p_2 \cdot(S'' p_2) - \overline{U}' \cdot S' + \frac12 V' \cdot V' + G^{(3)}_t \right] + O(\eps^3),
\]
we choose
\[
	G^{(3)}(x_2,p_3,\tau) = \int \left[\overline{U}'(x_2) \cdot S'(x_2,\tau) - p_3 \cdot \big(S''(x_2,\tau) p_3\big) + \frac12 \overline{V' \cdot V'}(x_2) - \frac12 V' \cdot V'(x_2,\tau)  \right] d\tau
\]
with the constant of integration such that $\overline{G^{(3)}} = 0$. Letting $A(x_2,\tau) = \int S(x_2,\tau) d\tau$, with $\bar{A} = 0$, we have
\[
	G^{(3)}(x_2,p_3,\tau) = \overline{U}'(x_2) \cdot A'(x_2,\tau) - p_3 \cdot \big(A''(x_2,\tau) p_3\big)
	+ \frac12 \int \left[ \overline{V' \cdot V'}(x_2) - (V' \cdot V')(x_2,\tau)\right] d\tau.
\]

Then the transformation   (\ref{eq:3})   takes form
\begin{align}
\begin{split}
	x_3 &= x_2 - 2\eps^3 A''(x_2,t/\eps) p_3 \\
	p_2 &= p_3 + \eps^3 G^{(3)}_x(x_2,p_3,t/\eps)
\end{split}
\end{align}
and the corresponding Hamiltonian is
\begin{align}
	H_3 = \frac12 p_3^2 + \overline{U}(x_3) + \frac{\eps^2}{2} \overline{V' \cdot V'}(x_3) + \eps^3 \left[ G^{(3)}_x \cdot p_3 - \overline{U}'(x_3) \cdot G^{(3)}_p - (S'' V')\cdot p_3 \right] + O(\eps^4)
	\label{eq:H3}
\end{align}
where we have used
\begin{align*}
	x_2 &= x_3 - \eps^3 G^{(3)}_p(x_3,p_3,t/\eps) + O(\eps^6) \\
	p_2 &= p_3 + \eps^3 G^{(3)}_x(x_3,p_3,t/\eps) + O(\eps^6),
\end{align*}
obtained from   (\ref{eq:3}). 
Finally, to eliminate the time dependence in the bracketed term in \eqref{eq:H3} we make one more transformation  
\begin{align}
	\begin{split}
	x_4 &= x_3 + \eps^4 G^{(4)}_p(x_3,p_4,t/\eps) \\
	p_3 &= p_4 + \eps^4 G^{(4)}_x(x_3,p_4,t/\eps)  
	\end{split}
	\label{eq:34}
\end{align}
 yielding the new Hamiltonian 
 \[
	H_4=H_3+ \varepsilon ^3 G^{(4)}_\tau.
\]  
As before, one can choose $ G^{(4)} $ to  replace the bracketed term in   (\ref{eq:H3})  by its average; since $ \overline{G^{(3)}}=0 $,  only the last terms in brackets survives, and 
\begin{align}
	H_4(x_4,p_4,t/\eps) = \frac12 p_4^2 + \overline{U}(x_4) + \frac{\eps^2}{2} \overline{V' \cdot V'} - \eps^3 \big(\overline{S'' V'}\big) \cdot p_4 + O(\eps^4).
	\label{eq:H4}
\end{align}
Combining the transformations, we obtain
\begin{align}\label{x3trans}
	x_4 = x_2 - 2\eps^3 A''(x_2,\tau) p_3 = x + \eps^2 S'(x,\tau) - 2\eps^3 A''(x,\tau) \dot x + O( \varepsilon^4);
\end{align}
we do not provide an explicit form for the coefficient of $ \varepsilon ^4 $ since it is rather cumbersome. 
This completes the proof of Theorem \ref{thm:1}
 \hfill $\diamondsuit$

\begin{rem}\label{rem:linearcase}
For the rotating quadratic saddle described in Section \ref{sec:examples}, the $O(\eps^3)$ term $G_x \cdot p_3 - \overline{U}'(x_3) \cdot G_p -(S'' V') \cdot p_3$ is time independent, hence does not need to be averaged. In other words, the coordinate transformation does not contain an $\eps^4$-term, and so $x_4 = x_3$.
By contrast, in the general nonlinear case it is  necessary to average the $\eps^3$-term and to compute the corresponding transformation for $x_4$.
\end{rem}

\subsection{Proof of Theorem \ref{thm:ode}}   
 Hamilton's equations associated with the Hamiltonian (\ref{eq:avhamiltonian})  are 
\begin{align}
	\begin{split}
	\dot X&=   P - \varepsilon ^3 \big(\overline{S'' V'} \cdot P\big)_P +O( \varepsilon ^4 )    \\
	\dot P &= -U^\prime - \varepsilon^2 W ^\prime + \eps^3 \big(\overline{S'' V'} \cdot P\big)^\prime +O( \varepsilon ^4 )  ,   
	\end{split}
	\label{eq:XP}
\end{align}
where $ (\cdot)_P $ denotes $P$-derivative of $ (\cdot) $. 
It is tempting to differentiate the first of the above equations and substitute into the result the expression for $ \dot P $ from the second. The problem, however, is that the remainder in the first equation in   (\ref{eq:XP})  depends on $ t/ \varepsilon $ and differentiation lowers the  order of magnitude to $ \varepsilon ^3 $. We remedy this problem by making one more transformation of the Hamiltonian \eqref{eq:avhamiltonian} via
\begin{align}
\begin{split}
	X_1 &= X + \eps^5 G^{(5)}_p(X,P_1,t/\eps)  \\
	P &= P_1 + \eps^5 G^{(5)}_x(X,P_1,t/\eps) , 
\end{split}
	\label{eq:XP1}
\end{align}
choosing $G^{(5)}$ so as to kill the time dependence in the coefficient of  $   \varepsilon ^4    $  in the  remainder of $H$.    Now   (\ref{eq:XP})  
transforms into 
\begin{align}
	\begin{split}
	\dot X&=   P - \varepsilon ^3 \big(\overline{S'' V'} \cdot P \big)_{P} +\varepsilon ^4 E(X,P)+O( \varepsilon ^5 )    \\
	\dot P&=  -W ^\prime + \eps^3 \big(\overline{S'' V'} \cdot P \big)^\prime +\varepsilon ^4 F(X,P)+O( \varepsilon ^5 )  ,   
	\end{split}
	\label{eq:XP5}
\end{align}
where we have dropped the subscripts, writing $ X, \,P $ instead of $ X_1,\, P_1 $ and where $E$ and $F$ are time independent. 

Observe that $\big(\overline{S'' V'} \cdot P \big)_P = \overline{S'' V'}$.  Differentiating the first equation in   (\ref{eq:XP5}) now yields
\begin{equation} 
	\ddot X =  \dot  P - \varepsilon ^3  (\overline{S'' V'} )^\prime \dot X+O( \varepsilon ^4 ). 
	\label{eq:ddoX}
\end{equation}  
Note that $\big(\overline{S'' V'} \cdot P \big)^\prime = \big((\overline{S'' V'})^\prime \big)^TP $, where $ (\cdot)^T $ denotes the transpose.  Substitution of the second equation in \eqref{eq:XP} into \eqref{eq:ddoX} leads to the claim  (\ref{eq:averagedequation}) and completes the proof of Theorem \ref{thm:ode}. \hfill $\diamondsuit$

\end{document}